
\magnification=1200
\def\div{{\rm div}}
\def\tr{{\rm tr}}

\def\f#1#2{{\textstyle{#1\over #2}}}

\def\next{\hfil\break\noindent}

\font\title=cmbx12

\noindent
{\title Compact null hypersurfaces and collapsing Riemannian manifolds}

\vskip .5cm\noindent
Alan D. Rendall\footnote*
{Present address: Max-Planck-Institut f\"ur
Gravitationsphysik, Schlaatzweg 1, 14473 Potsdam, Germany}
\next
Institut des Hautes Etudes Scientifiques
\next
35 Route de Chartres
\next
91440 Bures sur Yvette
\next
France
\next

\vskip .5cm\noindent
{\bf Abstract}

\noindent
Restrictions are obtained on the topology of a
compact divergence-free null hypersurface in a four-dimensional
Lorentzian manifold whose Ricci tensor is zero or satisfies some weaker
conditions. This is done by showing that each null hypersurface of
this type can be used to construct a family of three-dimensional
Riemannian metrics which collapses with bounded curvature and
applying known results on the topology of manifolds which collapse.
The result is then applied to general relativity, where it implies
a restriction on the topology of smooth compact Cauchy horizons in
spacetimes with various types of reasonable matter content.

\vskip .5cm

\noindent
{\bf 1. Introduction}

The concept of collapsing Riemannian manifolds introduced by
Cheeger and Gromov [3, 4] has been the subject of a considerable
amount of work in Riemannian geometry in recent years. A sequence
of Riemannian manifolds is said to collapse if the injectivity
radius tends uniformly to zero while the sectional curvature
remains bounded. A manifold is said to collapse if it admits a
collapsing sequence of Riemannian metrics. If this sequence can
be chosen so that the diameter of the Riemannian metrics remains
bounded then the manifold is said to collapse with bounded diameter.
In this paper results on collapsing Riemannian manifolds are applied
to a question in Lorentzian geometry which is of interest in
general relativity. A Lorentzian manifold $(M,g)$ is defined to be
a four-dimensional manifold $M$ with a metric tensor $g$ of
signature $(-,+,+,+)$. Hyperplanes in the tangent space to $M$
at some point $p$ are of three types; a hyperplane is said to
be spacelike, null or timelike if the pull-back of the metric
$g(p)$ to this hyperplane is positive definite, degenerate or
indefinite, respectively. A hypersurface is said to be spacelike,
null or timelike if its tangent space at each point is of the
corresponding type. If $H$ is a null hypersurface in $M$ and
if $H$ can be approximated by spacelike hypersurfaces in some
appropriate sense (the precise sense which is of relevance in
the following will be specified later) then a one-parameter family
of Riemannian metrics is obtained which degenerates as the
parameter tends to some limiting value. However this is not
enough to say that this family of metrics collapses; it is
also necessary to know that the curvature of these metrics is
uniformly bounded. In the following a situation will be exhibited where
this condition is satisfied.

Let $H$ be a smooth null hypersurface in a Lorentzian manifold. If $p$
is a point of $H$ there is precisely one direction in the tangent
space to $H$ at $p$ such that if a vector $L$ spans this direction
then $g(L,X)=0$ for all vectors $X$ tangent to $H$ at $p$. In fact
this is just the normal direction to $H$.  The curves in $H$ which are
tangent to this preferred direction are called the generators of $H$.
They are geodesics (see e.g. [20], p. 65).  Locally on $H$ it is
possible to choose a smooth vector field $L$ which on each generator
coincides with the tangent vector to an affinely parametrized
geodesic. This vector $L$ is only defined up to a constant
multiplicative factor on each generator. Let $\tilde L$ be a smooth
extension of $L$ to a null vector defined on an open neighbourhood in
$M$ of the domain of $L$ and define $\theta$ to be the restriction of
$\div \tilde L$ to $H$.  The quantity $\theta$ is called the expansion
of $L$ and does not depend on the extension $\tilde L$ chosen. If $L$
is scaled by a factor which is constant on each generator then
$\theta$ is scaled by the same factor. A null hypersurface is called
divergence-free if $\theta=0$ for some (and hence any) choice of $L$
in a neighbourhood of each point of $H$.

The main result of this paper is that in the presence of certain
restrictions on the Ricci tensor $r$ of the metric $g$, the topology
of a divergence-free compact null hypersurface in a Lorentzian
manifold is strongly constrained.

\noindent
{\bf Theorem 1} Let $(M,g)$ be a Lorentzian manifold whose Ricci tensor
is zero and let $H$ be a divergence-free compact smooth null hypersurface
in $M$. Then the manifold $H$ collapses with bounded diameter.

\vskip 10pt\noindent
This is a special case of the following more general result.

\noindent
{\bf Theorem 2} Let $(M,g)$ be a Lorentzian manifold and let $H$ be
a divergence-free compact smooth null hypersurface in $M$. Suppose
that $r(Z,X)=0$ when $Z$ is normal to $H$ and $X$ is tangent to $H$.
Then the manifold $H$ collapses with bounded diameter.

\vskip 10pt\noindent
To understand the significance of the conclusion it is necessary
to know how much the topology of a three-dimensional manifold
is constrained by the condition that it collapses. An introductory
account of the topological restrictions implied by collapsing can
be found in [14]. An example of a class of compact three-dimensional
manifolds which do not collapse are those which are hyperbolic, i.e.
which admit a metric of constant negative curvature. It turns out
that collapse with bounded diameter is a much stronger restriction
that just collapse. A compact three-dimensional manifold which
collapses with bounded diameter must either be a Seifert manifold
or admit a geometric structure of type Sol. Recall that a Seifert
manifold is one which admits a foliation by circles [18]. Any
such manifold admits a geometric structure in the sense of Thurston.
More precisely, the class of Seifert manifolds coincides with the
class of manifolds which admit a geometric structure corresponding
to one of six of the eight Thurston geometries. The two remaining
geometries are Sol, which was mentioned above, and hyperbolic geometry,
which is not compatible with collapsing. Thus a compact three-dimensional
manifold which collapses with bounded diameter admits a geometric
structure which is not hyperbolic. All other geometric structures
can occur, as discussed in more detail in section 3.

The argument which leads to the statement on the existence of a
geometric structure will now be sketched. (Cf. [16], [19].) This uses
the notion of Gromov-Hausdorff convergence of metric spaces (see e.g.
[14]). The metric spaces of relevance here are obtained by considering
the distance function associated to a Riemannian metric. The Gromov
precompactness theorem [10] says that a sequence of Riemannian metrics
with bounded curvature and diameter has a subsequence which converges
in the Gromov-Hausdorff sense to some metric space. A theorem of
Fukaya ([7], Theorem 0.6) then shows that the limiting space is
isometric to the quotient of a Riemannian manifold by an isometric
action of $O(n)$. This means in particular that the limiting space
will have an open dense subset which is a Riemannian manifold. This
open set corresponds to the principal orbit type of the $O(n)$ action
[13]. For the present purposes the dimension of the limiting space
will be defined to be the dimension of this open dense subset. There
are now three possibilities: the dimension of the limiting space may
be zero, one or two. If the dimension is zero then $H$ is by
definition an almost flat manifold [9, 1, 17].  Any almost flat
three-dimensional manifold either admits a flat metric or a geometric
structure of type Nil. If the limiting space is of dimension two then
by a theorem of Fukaya ([6], Proposition 11.5) it is an orbifold. Then
another result of Fukaya (the fibre bundle theorem) [8] shows that $H$
is a Seifert manifold. Finally, if the limiting space has dimension
one, it is homeomorphic either to a circle or a closed interval. In
particular it is an orbifold.  A result of Tuschmann [19] then implies
that it must be flat or admit a geometric structure of type Nil or
Sol.

\vskip 10pt
Theorems 1 and 2 are proved in section 2. The idea is to approximate
the null hypersurface $H$ by spacelike hypersurfaces in an
appropriate sense and show that the family of Riemannian metrics
which arises in this way has bounded curvature. This can then
be used to show that $H$ collapses. In section 3 the applications of the
theorem to general relativity are described. In general relativity there
is an important class of compact null hypersurfaces, the compact Cauchy
horizons, which are divergence-free. The hypotheses on the curvature in
Theorems 1 and 2 follow from certain hypotheses on the matter content of
spacetime. In particular, Theorem 1 corresponds to the case where
spacetime is empty.

\vskip 10pt\noindent
{\bf 2. Boundedness of the curvature}

Let $H$ be a smooth compact null hypersurface in a Lorentzian
manifold $(M,g)$. There exists a one-dimensional distribution
(field of directions) on $M$ which is timelike in the sense
that it can be spanned locally by a vector $T$ with $g(T,T)<0$.
Let $W$ be an open neighbourhood of $H$ with compact closure. If
$g$ is time orientable then there exists a global smooth vector
field $T$ on $W$ spanning the given distribution. It may be assumed
without loss of generality that $g(T,T)=-1$. In this case
a family of tensors is defined by:
$$g_\lambda (X,Y)=g(X,Y)+\lambda g(T,X)g(T,Y)\eqno(1)$$
For $\lambda$ sufficiently close to zero $g_\lambda$ is Lorentz
metric. Even if $g$ is not time orientable it is still possible
to do a similar construction. The choice of the unit vector
$T$ at a given point is unique up to a sign and the definition
of $g_\lambda$ does not change when $T$ is replaced by $-T$.
Hence if $g_\lambda$ is defined by (1) using local vector fields $T$,
these locally defined tensors fit together to define a smooth
object $g_\lambda$ globally on $W$. If $X$ is a non-zero vector tangent
to $H$ then $g(X,X)\ge 0$. It follows that $g_\lambda(X,X)>0$ for
$\lambda>0$. Thus $H$ is spacelike with respect to the Lorentz metrics
$g_\lambda$ for each $\lambda >0$. Let $h_\lambda$ be the pull-back of
$g_\lambda$ to $H$. Then $h_\lambda$ is a Riemannian metric for
$\lambda>0$ and is degenerate for $\lambda=0$. The diameter of
the metric $h_\lambda$ is bounded as $\lambda\to 0$ while its
volume tends to zero.

It will now be shown that the curvature of $h_\lambda$ remains
bounded as $\lambda\to 0$. The curvature will be computed in a
local frame adapted to $H$. This frame will only be defined at
points of $H$. Let $Z$ be a vector field defined on
a neighbourhood of a point of $H$ which is non-vanishing and tangent
to the generators of $H$. Consider a unit timelike vector field
$T$ as above. Let $\Pi$ denote the two-dimensional distribution
which is the intersection of the tangent space of $H$ with the normal
to $T$ with respect to $g_\lambda$. It follows from the definition
of $g_\lambda$ that $\Pi$ does not depend on $\lambda$ and that
the induced metric on $\Pi$ is also independent of $\lambda$. Let
$X$ and $Y$ be an orthonormal basis of $\Pi$. Let $z=g_\lambda (Z,Z)$.
This is positive for $\lambda>0$ and tends to zero as $\lambda\to 0$.
For $\lambda>0$ let $\hat Z=z^{-1/2}Z$ and let $\hat U$ be the unit
normal vector to $H$. It can be rescaled to give a vector $U$ which
extends smoothly to $\lambda=0$. The vector $U$ is proportional to
$Z$ for $\lambda=0$ and it can be assumed without loss of generality
that $U=Z$ there. Note that for $\lambda=0$:
$$dz/d\lambda=(g(T,Z))^2>0$$
It follows that $U-Z=zD$ for some smooth vector $D$. Let
$u=-g_\lambda(U,U)$. Then $u=z-z^2g_\lambda(D,D)$. In particular,
it follows that $u/z$ has a smooth extension to $\lambda=0$ and
takes the value one there.
Note that since the generators are geodesics
there is a smooth function $f$ and a smooth vector field $W$ such that
$\nabla_Z Z=fZ+zW$.

Let $k_\lambda$ denote the second fundamental form of the hypersurface
$H$ with respect to the metric $g_\lambda$. Let $P$ and $Q$ be any
vectors in the plane spanned by $X$ and $Y$.
$$\eqalign{
k_\lambda(P,Q)&=(z/u)^{1/2}
[z^{-1/2}g_\lambda(\nabla_P Q,Z)+z^{1/2}g_\lambda(\nabla_P Q, D)]   \cr
k_\lambda(\hat Z,P)&=(z/u)^{1/2}
[-g_\lambda (W,P)+g_\lambda(\nabla_Z P,D)]     \cr
k_\lambda(\hat Z,\hat Z)&=u^{-1/2}g_\lambda(W,U)}$$
The curvature of $h_\lambda$ can be expressed in terms of the curvature
of $g_\lambda$ and $k_\lambda$ using the Gauss equation. Many
of the components of $k_\lambda$ appear to blow up as $\lambda\to 0$.
We will see that in the case of a compact null hypersurface several
apparently divergent terms vanish.

\noindent
{\bf Lemma 1} Under the hypothesis of Theorem 2 the family of metrics
$h_\lambda$ has bounded curvature.

\noindent
{\bf Proof} Note first that since the dimension is
three, in order to bound the sectional curvature of $h_\lambda$
it suffices to show the boundedness of the components of its
Ricci tensor $p$ in an orthonormal frame. The Gauss
equation gives:
$$\eqalign{
p(P,Q)&=r(P,Q)-(\tr k)k(P,Q)+(k\cdot k)(P,Q)+u^{-1}g(R(U,Q)U,P))
            \cr
p(\hat Z,P)&=z^{-1/2}r(Z,P)-(\tr k)k(\hat Z,P)+(k\cdot k)(\hat Z,P)
+u^{-1}z^{-1/2}g(R(U,Z)U,P))          \cr
p(\hat Z,\hat Z)&=z^{-1}r(Z,Z)-(\tr k)k(\hat Z,\hat Z)
+(k\cdot k)(\hat Z,\hat Z)+u^{-1}z^{-1}g(R(U,Z)U,Z)
}$$
The aim is to show that the quantities $p(P,Q)$, $p(\hat Z,P)$ and
$p(\hat Z,\hat Z)$ are bounded as $\lambda\to 0$. Note first that
under the hypotheses of Theorem 2 the contributions of $r$ in the
above expressions are bounded. Consider next the contributions of the
spacetime curvature tensor.
$$\eqalign{g(R(U,Z)U,P)&=-zg(R(U,D)U,P)       \cr
g(R(U,Z)U,Z)&=z^2g(R(U,D)U,D)}$$
Thus all contributions of the curvature tensor will be bounded if
we can show that $z^{-1}g(R(U,\cdot)U,\cdot)$ gives something bounded
when evaluated on any pair of regular vectors. To do this it is
necessary to study the behaviour of the null geodesics which are
tangent to $Z$ for $\lambda=0$. For this computation the vector
$Z$ will be chosen to be parallelly transported along the null geodesics
lying in $H$ for $\lambda=0$. The expansion corresponding to $Z$ is
given by $\theta= g(\nabla_X Z,X)+g(\nabla_Y Z,Y)$. Define a tensor
$B$ by $B(F,G)=g(\nabla_F Z, G)$. Then, for $\lambda=0$, the trace of
$B$ is equal to $\theta$ and $B(F,G)=B(G,F)$ whenever $F$ and $G$ are
tangent to $H$. The evolution of $\theta$ is given by the
Raychaudhuri equation ([20], p. 222):
$$\nabla_Z\theta=-[B(X,X)]^2-[B(Y,Y)]^2-2[B(X,Y)]^2-r(Z,Z)$$
By hypothesis $\theta=0$ and $r(Z,Z)=0$. It follows that
$B(X,X)=B(X,Y)=B(Y,Y)=0$. Combining this with the condition that
$Z$ is parallelly transported shows that $B(F,G)=0$
for any $F$, $G$ tangent to $H$. Now suppose that $F$ and $G$ are
vectors tangent to $H$ which are parallelly transported along the
integral curves of $Z$. The evolution of $B(F,G)$ is described by
the equation
$$\nabla_Z(B(F,G))=-B(F,X)B(G,X)-B(F,Y)B(G,Y)+g(R(Z,F)Z,G)$$
It follows that $g(R(Z,F)Z,G)=0$ and hence $g(R(U,F)U,G)$ vanishes
for $\lambda=0$. This bounds the contributions of the curvature
tensor in the expressions for $p(F,G)$. It remains to consider
the contributions of the second fundamental form. Note that
the discussion of the Raychaudhuri equation allows us to conclude
and that $k(P,Q)$ is $O(z^{1/2})$ as $\lambda\to 0$.
Also $\tr k$ is $O(z^{-1/2})$. It follows
immediately that $p(P,Q)$ is bounded. To see that $p(\hat Z,P)$
is bounded it suffices to observe that both $(\tr k)k(\hat Z,P)$
and $(k\cdot k)(\hat Z, P)$ differ from $k(\hat Z,\hat Z)k(\hat Z, P)$
by bounded quantities so that the singular terms cancel. Finally,
$$\eqalign{&-(\tr k)k(\hat Z,\hat Z)+(k\cdot k)(\hat Z,\hat Z) \cr
&=(-k(\hat Z,\hat Z)-k(X,X)-k(Y,Y))k(\hat Z,\hat Z)
+(k(X,\hat Z))^2+(k(Y,\hat Z))^2+(k(\hat Z,\hat Z))^2          \cr
&=-(k(X,X)+k(Y,Y))k(\hat Z,\hat Z)+(k(X,\hat Z))^2+(k(Y,\hat Z)^2}$$
so that $p(\hat Z,\hat Z)$ is also bounded.

\vskip 10pt\noindent
{\bf Proof of Theorem 2}
If a family of Riemannian metrics is such that the curvature is bounded
and the volume goes to zero then it follows from Bishop's theorem
[2] that the injectivity radius goes uniformly to zero. Hence it can be
concluded from Lemma 1 that the family of metrics $h_\lambda$ associated
to a divergence-free compact null hypersurface as above does represent a
collapse in the sense of Cheeger and Gromov and that the manifold $H$
collapses.

\vskip 10pt\noindent
{\bf 3. Compact Cauchy horizons}

In this section theorem 2 will be applied to general relativity.
A Lorentz manifold $(M,g)$ satifies the Einstein equations if
$r=T-\f12\tr T g$, where $T$ is the energy-momentum tensor. What
exactly this energy-momentum tensor is depends on the assumed matter
content of spacetime. The tensor $T$ is said to satisfy the weak
energy condition if $T(V,V)\ge 0$ for any timelike vector $V$.
The following result can be deduced from theorem 2. For definitions
of the concepts used in its statement, see [11].

\noindent
{\bf Theorem 3} Let $(M,g)$ be a solution of the Einstein equations
with energy-momentum tensor $T$ containing a partial Cauchy surface
with a smooth compact Cauchy horizon $H$. Suppose that:
\next
(i) $T$ satifies the the weak energy condition
\next
(ii) if $T(N,N)=0$ for a null vector $N$ then $T(N,X)=0$ for any
vector $X$ orthogonal to $N$
\next
Then $H$ collapses with bounded diameter.

\noindent
{\bf Proof} Since $H$ is a smooth compact Cauchy horizon it has
zero divergence ([11], pp. 295-298). Applying the Raychaudhuri
equation and using condition (i) above shows that $r(Z,Z)=0$.
Hence, by condition (ii) $r(Z,X)=0$ for any $X$. Thus the
hypotheses of theorem 2 are satisfied.

\vskip 10pt\noindent
In order that this theorem be interesting it is necessary to
check that it is satisfied for types of energy-momentum tensor
which are important in general relativity. Some examples will
now be given.

\noindent
{\bf Example 1} $T=0$. This is the case of the vacuum Einstein
equations.
\next
{\bf Example 2} $T=(\rho+p)U\otimes U+pg$, where $U$ is a unit timelike
vector and $\rho$ and $p$ are positive functions
with $p<\rho$. This is the case of a perfect fluid where the speed
of sound is less than the speed of light.
\next
{\bf Example 3} $T$ is the energy-momentum tensor of a kinetic
matter model. See [15] for further discussion of this and
the other matter models above.
\next
{\bf Example 4} $T=F\cdot F-\f14 (|F|^2)g$, where $F$ is a
two-form. This is the energy-momentum tensor of an electromagnetic
field.

Theorems 1 and 2 give necessary conditions for a compact manifold
to occur as a compact Cauchy horizon. This will now be complemented
by mentioning examples of manifolds which do occur as Cauchy horizons
in solutions of the vacuum Einstein equations. It fact it follows
from [5] that manifolds occur which admit a geometric structure
belonging to any Thurston geometry except Sol and hyperbolic. Since
hyperbolic manifolds have been ruled out, the one geometry for which
neither positive or negative results are known is Sol. Manifolds
with a geometric structure of type Sol can collapse [6]. What is not
clear is whether this collapse can be realized by a Cauchy horizon
in a solution of the Einstein equations. This issue is related to
a conjecture of Isenberg and Moncrief [12]. This says that any analytic
vacuum spacetime containing a compact Cauchy horizon must admit
at least one Killing vector with closed orbits which is tangent
to the horizon (possibly under some additional hypotheses).
This means in particular that $H$ admits a foliation
by circles, i.e. that it is a Seifert manifold [18]. This in turn
means that it admits a geometric structure which is not Sol or
hyperbolic. In other words, if the Isenberg-Moncrief conjecture
is correct then a manifold with a geometric structure of type
Sol cannot occur as a compact Cauchy horizon (at least under the
assumptions of analyticity and the vacuum Einstein equations). Thus
Theorems 1 and 2 would not capture all restrictions on the topology of
a compact Cauchy horizon. On the other hand, pushing further the idea
of applying collapsing to Cauchy horizons might be helpful in
proving the conjecture, which has as yet only been proved in a
special case.

\vskip 10pt\noindent
{\bf Acknowledgements} I thank Dennis Sullivan for pointing out that
collapsing might be useful in the study of Cauchy horizons. I am
also grateful to John Lott, Xiaochun Rong and Wilderich Tuschmann for
helpful correspondence and discussions.

\vskip 10pt\noindent
{\bf References}

\noindent
[1] Buser, P., Karcher, H.: Gromov's almost flat manifolds. Ast\'erisque
81, 1-148 (1981).
\next
[2] Chavel, I.: Riemannian geometry: a modern introduction. (Cambridge
University Press, 1993).
\next
[3] Cheeger, J., Gromov, M.: Collapsing Riemannian manifolds while
keeping their curvature bounded I. J. Diff. Geom. 23, 309-346 (1986).
\next
[4] Cheeger, J., Gromov, M.: Collapsing Riemannian manifolds while
keeping their curvature bounded II. J. Diff. Geom. 32, 269-298 (1990).
\next
[5] Chru\'sciel, P., Rendall, A. D.: Strong cosmic censorship in
vacuum spacetimes with compact locally homogeneous Cauchy surfaces.
To appear in Ann. Phys. (NY)
\next
[6] Fukaya, K.: Hausdorff convergence of Riemannian manifolds and
applications. In T. Ochiai (ed.) Recent topics in differential and
analytic geometry. Academic Press, Boston, 1990.
\next
[7] Fukaya, K.: A boundary of the set of Riemannian manifolds with bounded
curvatures and diameters. J. Diff. Geom. 28, 1-21 (1988).
\next
[8] Fukaya, K.: Collapsing Riemannian manifolds to ones of lower dimension II.
J. Math. Soc. Japan 41, 333-356 (1989).
\next
[9] Gromov, M.: Almost flat manifolds. J. Diff. Geom. 13, 231-241 (1978).
\next
[10] Gromov, M. (r\'edig\'e par J. Lafontaine et P. Pansu) Structures
m\'etriques pour les vari\'et\'es riemanniennes. Cedic/ Fernand
Nathan, Paris, 1981.
\next
[11] Hawking, S. W., Ellis, G. F. R.: The large-scale structure of
space-time. Cambridge University Press, 1973.
\next
[12] Isenberg, J., Moncrief, V.: On spacetimes containing Killing vector
fields with non-closed orbits. Class. Quantum Grav. 9, 1683-1691 (1992).
\next
[13] J\"anich, K.: Differenzierbare G-Mannigfaltigkeiten. Lecture
Notes in Mathematics 59 (Springer, Berlin, 1968).
\next
[14] Pansu, P.: Effondrement des vari\'et\'es riemanniennes (d'apr\`es
J. Cheeger et M. Gromov). Ast\'erisque 121-122, 63-82 (1985).
\next
[15] Rendall, A. D.: Global properties of locally homogeneous
cosmological models with matter. To appear in Math. Proc. Camb. Phil. Soc.
\next
[16] Rong, X.: The limiting $\eta$-invariants of collapsed 3-manifolds.
J. Diff. Geom. 37, 535-568 (1993).
\next
[17] Ruh, E. A.: Almost flat manifolds. J. Diff. Geom. 17, 1-14 (1982).
\next
[18] Scott, P.: The geometries of 3-manifolds. Bull. London Math.
Soc. 15, 401-487 (1983).
\next
[19] Tuschmann, W.: Collapsing, solvmanifolds and infrahomogeneous
spaces. Preprint IHES/M/95/7.
\next
[20] Wald, R.: General relativity. (University of Chicago Press,
Chicago, 1984).

\end